# On a new measure on infinite dimensional unite cube

## I.Sh. Jabbarov

Ganja State University (Azerbaijan)

*e-mail: ilgar_j@rambler.ru*

**Abstract.** Measure Theory plays an important role in many questions of the Mathematics. The notion of a measure being introduced as a generalization of a notion of a segment's size made many of limiting processes be a formal procedure, and by this reason stood very productive in the questions of Harmonic analysis ([1-6]). Discovery of the Haar measure was a valuable event for the harmonic analysis in the topological groups. It stood clear that many of measures, particularly, the product of Lebesgue measure could be considered as a Haar measure. The product measure has many important properties concerning projections (see [1,3]). The theorems of Fubini and Tonelly made it very useful in applications.

In this work we introduce a new measure in infinite dimensional unite cube $[0,1] \times [0,1] \times \cdots$. We use the Tichonoff metric to define a set function in the algebra generated by open balls defining their measure in product meaning. By this way we introduce a new measure in infinite dimensional unite cube different from the Haar or product measures and discuss some its properties.

Main difference between the introduced measure and product measure connected with the property: *let we are given with a infinite family of open balls every of which does not contain any other with total finite measure; then they have an empty intersection*. Consequently, every point contained in by a finite number of considered balls only.

This property does not satisfied by cylinder sets. For example, let
$$D_1 = I_1 \times I \times I \times \cdots, \quad D_2 = I_2 \times I_1 \times I \times \cdots, \ldots$$
$$I = [0,1], \quad I_k = [0, k^2/(k+1)^2), \quad k = 1, 2, \ldots$$

It clear that every of these cylinder sets does not contain any other, but their intersection is not empty (contains the point $(0,0,\ldots)$). This makes two measures currently different.

## 1. Introduction.

Let's enter the Tichonoff metric in the unite cube $\Omega = \{(x_n) \mid 0 \leq x_n \leq 1, n = 1, 2, \ldots\}$ as below:

$$d(x.y) = \sum_{n=1}^{\infty} e^{1-n} |x_n - y_n|; \qquad (1)$$

here $x = (x_n), y = (y_n) \in \Omega$. Let's define the ball of a radius $r > 0$ and a centre $\theta \in \Omega$ by the equality

$$B(\theta, r) = \{x \in \Omega \mid d(x, \theta) < r\}.$$

It is best known that in the cube $\Omega$ a product Lebesgue measure may be introduced (see [1, p. 219]). There is also another construction of a measure called a Haar measure. The Haar measure is a measure defined in the locally compact topological groups. It was proven also uniqueness of this measure (see [2, p.241]). Many of measures used in various brunches of the mathematics could be considered as a Haar measure. Particularly, the product of Lebesgue measures in $\Omega$ is a Haar measure (with an errow of null set).

Our goal here is to define a *new measure* different from mentioned above measures. For the construction we shall use open balls instead of cylinder sets. Despite that they generate the same $\sigma$-algebra we get a different measure. In the literature it was observed that measurable in Haar meaning sets, in general, does not allow appoximations by open balls (see [7]). The constructed example is very special. Here we show that the case not casual, and show the cause way this happens in infinite dimensional unite cube (see lemma 2 below).

## 2. Construction of a new measure

**Definition 1.** *Let $\sigma : N \to N$ be any one to one mapping of the set of natural numbers. If for any $n > m$ there is a natural number $m$ such that $\sigma(n) = n$, then we call $\sigma$ a finite permutation. A subset $A \subset \Omega$ is called to be finite-symmetric if for any element $\theta = (\theta_n) \in A$ and any finite permutation $\sigma$ one has $\sigma\theta = (\theta_{\sigma(n)}) \in A$.*

Let $\Sigma$ to denote the set of all finite permutations. We shall define on this set a product of two finite permutations as a composition of mappings. Then $\Sigma$ becomes a group which contains each group of $n$ degree permutations as a subgroup (we consider each $n$ degree permutation $\sigma$ as a finite permutation in the sense of definition 1, i. e. we put $\sigma(m) = m$ when $m > n$). The set $\Sigma$ is a countable set and we can arrange its elements in a sequence.

Let $\omega \in \Omega$, $\Sigma(\omega) = \{\sigma\omega \mid \sigma \in \Sigma\}$ and $\Sigma'(\omega)$ to mean the closed set of all limit points of the sequence $\Sigma(\omega)$. For real $t$ we denote $\{t\Lambda\} = (\{t\lambda_n\})$ where $\Lambda = (\lambda_n)$. Let $\mu$ to denote the product of linear Lebesgue measures $m$ given on the interval $[0,1]$: $\mu = m \times m \times \cdots$.

To construct the measure let's consider the open ball
$$B(0,r) = \{x \in \Omega_0 \mid d(x,0) < r\}$$
in the cube $\Omega_0 = \{x = (x_n) \mid |x_n| \leq 1\}$. Since $|x_n| \leq 1$ then for the natural number $N$ we have
$$\sum_{n=N+1}^{\infty} e^{1-n} |x_n| \leq e^{-N} \sum_{n=0}^{\infty} e^{-n} < e^{1-N}.$$
Taking arbitrarily small real number $\varepsilon > 0$ we get
$$\sum_{n=1}^{N} e^{1-n} |x_n| \leq d(x,0) \leq \sum_{n=1}^{N} e^{1-n} |x_n| + \varepsilon$$
when $N \geq \log e\varepsilon^{-1}$. Therefore,
$$B_N(0, r-\varepsilon) \times [0,1] \times \cdots \subset B(0,r) \subset B_N(0,r) \times [0,1] \times \cdots,$$

where $B_N(0,r)$ denotes the projection of the ball $B(0,r)$ to the subspace of first $N$ coordinates. Then, for the volume $\mu_N(r)$ of the projection $B_N(0,r)$ we have (see [8] or [9, p.319]):

$$\mu_N(r) - \mu_N(r-\varepsilon) = \int_{r-\varepsilon \leq \sum_{n=1}^{N} e^{1-n}|x_n| \leq r} dx_1 \cdots dx_N = 2^N \int_{r-\varepsilon \leq u \leq r} du \int_{\sum_{n=1}^{N} e^{1-n} u_n = u} \frac{ds}{\|\nabla\|} \leq$$

$$\leq \varepsilon 2^N \int_M \frac{ds}{\|\nabla\|},$$

and the last integral is an surface integral over the surface $M$ defined by the equation

$$\sum_{n=1}^{N} e^{1-n} u_n = u, \quad 0 \leq u_k \leq 1, 1 \leq k \leq N; \quad (2)$$

here $\nabla$ is a norm of the gradient of linear function on the left side of the latest equality, i.e.

$$\|\nabla\| = \sqrt{1 + e^{-2} + \cdots + e^{2-2N}}.$$

Defining $u_1$ from (2) we get

$$\int_M \frac{ds}{\|\nabla\|} \leq \int_0^1 \cdots \int_0^1 du_2 \cdots du_N = 1.$$

So, we have

$$\mu_N(r) - \mu_N(r-\varepsilon) \leq \varepsilon 2^N.$$

By taking the greatest $N$, satisfying the condition $N \geq \log e\varepsilon^{-1}$, i.e. $N = [\log e\varepsilon^{-1}] + 1$, we may write $\varepsilon \leq e^{2-N}$. Then from the last inequality it is follows that

$$\mu_N(r) - \mu_N(r-\varepsilon) \leq 2^N e^{2-N} \to 0,$$

as $N \to \infty$, or as $\varepsilon \to 0$. Since $B_{N+1}(0,r) \subset B_N(0,r) \times [0,1]$ then the sequence $(\mu_N(r))$ is monotonically decreasing. So, it is bounded below by $\mu_{N_0}(r/2)$, with $N_0 = [\log 2er^{-1}] + 1$. Therefore, there exists a limit

$$\lim_{\varepsilon \to 0} B_N(0, r-\varepsilon) = \lim_{N \to \infty} B_N(0,r) = \mu_0(B(0,r))$$

which we are accepting as a measure of the ball $B(0,r)$.

The measure of the ball in $\Omega$ we define as a limit of the measures of intersections $\Omega \cap B_N(\theta, r) \times [0,1] \times \cdots$ as $N \to \infty$, where $N$ defined for given $\varepsilon$ as above. The measure of the complement of the ball $B(\theta, r)$ is defined simply as $1 - \mu(B(\theta,r))$ or as a limit of measures of complements $[0,1]^N \setminus B_N(\theta, r-\varepsilon)$, as $\varepsilon \to 0$.

Consider now a union $A_1 \cup \cdots \cup A_k$ every component $A_i$ of which is some ball in $\Omega$ (naturally, with the different centers and radiuses) or it's complement. The measure of this union we define as above. Fixing $\varepsilon > 0$ we replace the components $A_i$ by cylinder sets with a "tower" $B_N(\theta, r)$, if it is a ball, or with a "tower"

$$[0,1]^N \setminus B_N(\theta, r - \varepsilon),$$

if it is a complement of the ball. We get then a union

$$\bigcup_i C_N(\theta_i, r_i),$$

with components $C_N(\theta_i, r_i)$ every of which being either $B_N(\theta, r)$ or $[0,1]^N \setminus B_N(\theta, r - \varepsilon)$. The set we have got has an $N$-dimensional Lebesgue measure. The $\varepsilon$ is taken less than $\min_i r_i$. As well as above, the error is estimated as a value $\leq k\varepsilon 2^N$. The demanded mesure we get after of passing to the limit.

Let's denonote by $\Pi$ the algebra generated by the class of open balls in $\Omega$. So, we have constructed some finite additive set function $\mu_0$ in $\Pi$. It is easy to note that $\mu_0$ is a restriction of a product measure into the $\Pi$. The algebra $\Pi$ can uniquely be extended to the $\sigma$-algebra $\Xi$ of subsets in $\Omega$ (see [1-6]). Corresponding continuation of a set function defined above to the $\sigma$-algebra $\Xi$ gives in $\Omega$ some measure and we get a measure space $(\Omega, \Xi, \mu_0)$. Clearly, the $\sigma$-algebra $\Xi$ coinsides with the the $\sigma$-algebra generated by cylinder sets.

**Definition 2.** A subset $A \subset \Omega$ is called to be a subset of zero measure if for any $\varepsilon > 0$ there exist open balls $B_1, B_2, \ldots$ satisfying conditions:

1) $A \subset B = B(\varepsilon) = \bigcup_{k=1}^{\infty} B_k$;

2) none of these balls contains any other;

3) $\sum_{k=1}^{\infty} \mu_0(B_k) < \varepsilon$.

The inner and outer measures $\mu_{0*}$ and $\mu_0^*$ are defined by a similar way by using of covering with open balls satisfying the natural condition 2) above. Since each open ball can be enclosed in some cylinder set with a measure enough close to the measure of the ball, we have for any subset $A \in \Omega$

$$\mu_{0*}(A) \leq \mu_*(A) \leq \mu^*(A) \leq \mu_0^*(A).$$

We call now a subset $A$ to be $\mu_0$-*measurable* in $\Omega$ if and only if the equality $\mu_{0*}(A) = \mu_0^*(A)$ is satisfied. Clearly, that the subset $A$ is measurable iff there exist a subset $D \subset \Xi$ for which the subset $D \Delta A$ has a zero measure. The defined Lebesgue extension of a measure $\mu_0$ is a regular measure and every $\mu_0$- measurable subset is measurable in the sense of product measure. Below we shall see that the definition of measureability differs from that in product meaning.

### 3. Supplementary basic results

**Lemma 1.** *Let A be a subset of zero measure in the meaning of the definition 2, $\varepsilon > 0$ given number, and open balls $B_1, B_2,...$ satisfy the conditions:*

*1) $A \subset B = \bigcup_{k=1}^{\infty} B_k$ ;*

*2) none of these balls contains any other;*

*3) $\sum_{k=1}^{\infty} \mu_0(B_k) < \varepsilon$ .*

*Further, let for given $a \in A$ we have*

$$\Sigma(a) \cup \Sigma'(a) \subset B.$$

*Then one can find a natural number m such that*

$$\Sigma(a) \cup \Sigma'(a) \subset \bigcup_{k=1}^{m} B_k,$$

*and*

$$(\Sigma(a) \cup \Sigma'(a)) \cap \left(\bigcup_{k=1}^{m} B_k\right) = \varnothing.$$

*Proof.* It is enough to show that one can find a natural number $m$ satisfying the relation

$$\Sigma(a) \subset \bigcup_{k=1}^{m} B_k .$$

Suppose the contrary. Let

$$\Sigma(a) \cap \left(\bigcup_{k=1}^{m} B_k\right) \neq \varnothing$$

for every natural $m$. Two cases are possible: 1) there exist a point $\bar{\theta} \in \Sigma(a)$ belonging to infinite number of balls; 2) there will be a subsequence of elements $\bar{\theta}_j, \theta_j \in \Sigma(a)$ which does not contained in by any finite union of balls $B_k$. We shall consider both possibilities separately and shall prove that they lead to the contradiction.

1) Let $B_{k_1}, B_{k_2}, B_{k_3},...$ be sequence of all balls into which the element $\bar{\theta}$ belongs. We shall denote $d$ the distance from $\bar{\theta}$ to the bound of $B_{k_1}$. As $B_{k_1}$ is open set, then $d > 0$. Let $B_k$ be any ball of a radius $< d/2$ from the list above, containing the point $\bar{\theta}$. From the told it follows that the ball $B_k$ should contained in the ball $B_{k_1}$. But it contradicts the agreement accepted above.

2) Let $\bar{\theta}$ be some limit point of the sequence $(\bar{\theta}_j)$. According to the condition of the lemma 3 $\bar{\theta} \in B_s$ for some $s$. Let $d$ denotes the distance from $\bar{\theta}$ to the bound of $B_s$. As $\bar{\theta}$ is a limit point, then a ball with the centre in the point $\bar{\theta}$ and radius $d/4$ contains an infinite set of members of the sequence $(\bar{\theta}_j)$, say members $\bar{\theta}_{j_1}, \bar{\theta}_{j_2}, \ldots$. According to 1), each point of this sequence can belong only to finite number of balls. So, the specified sequence will be contained in a union of infinite subfamily of balls $B_k$. Among them will be found infinitely many number of balls having radius $< d/4$. All of them, then, should contained in the ball $B_s$. The received contradiction excludes the case 2) also.

So, there exists such $n$ for which $\Sigma(a) \subset \bigcup_{k \leq n} B_k$. Lemma 1 is proved.

The lemma 1 is not satisfied by any denumerable family of cylider subsets, as it is clear from the example above with some modification (we take $a=(0, 0,\ldots)$).

**Lemma 2.** *Let $A \subset \Omega$ be a finite-symmetric subset of zero measure and $\Lambda = (\lambda_n)$ is an unbounded, monotonically increasing sequence of positive real numbers any finite subfamily of elements of which is linearly independent over the field of rational numbers. Let $B \supset A$ be any open, in the Tychonoff metric, subset with $\mu_0(B) < \varepsilon$,*

$$E_0 = \{0 \leq t \leq 1 \mid \{t\Lambda\} \in A \wedge \Sigma'\{t\Lambda\} \subset B\}.$$

*Then, we have $m(E_0) \leq c_0 \varepsilon$ where $c_0 > 0$ is an absolute constant, $m$ designates the Lebesgue measure.*

*Proof.* Let $\varepsilon$ be any small positive number. As the numbers $\lambda_n$ are linearly independent, for any finite permutation $\sigma$, one has $(\{t_1\lambda_n\}) \neq (\{t_2\lambda_{\sigma(n)}\})$ when $t_1 \neq t_2$; otherwise we should receive the equality $\{t_1\lambda_s\} = \{t_2\lambda_s\}$ for some natural $s$ which is invariant for $\sigma$. Then one has $(t_1 - t_2)\lambda_s = k, \ k \in Z$. Further, writing out the similar equality for natural $r > s$, we get the relation

$$k_1/\lambda_r - k/\lambda_s = \frac{k_1\lambda_s - k\lambda_r}{\lambda_r \lambda_s} = 0$$

for some integral $k_1$ which contradicts the linear independence of the numbers $\lambda_n$. Hence, for any pair of various numbers $t_1$ and $t_2$ one has $(\{t_1\lambda_n\}) \notin \{(\{t_2\lambda_{\sigma(n)}\}) \mid \sigma \in \Sigma\}$ (the outer brackets denote a set). By the conditions of the lemma 1, there exists a family of open balls $B_1, B_2,\ldots$ such that each ball does not contain any other one from this family (the ball contained in by other one can be deleted) and

$$A \subset B \subset \bigcup_{j=1}^{\infty} B_j, \sum \mu(B_j) < 1.5\varepsilon.$$

Now we take some permutation $\sigma \in \Sigma$ satisfying the equalities $\sigma(1) = n_1, \ldots, \sigma(k) = n_k$ where natural numbers $n_j$ are taken as below. Let $B'_N$ to denote the projection of the ball $B_1$, with $\mu(B_1) = \varepsilon_1$, into the subspace of first $N$ co-ordinate axes where the number $N$ is taken so that

$$\mu(B'_N) < 2\varepsilon_1.$$

Let $B'_N$ be enclosed by a union of cubes with edge $\delta$ and a total measure not exceeding $3\varepsilon_1$ having intersections over their boundary only. We put down $k = N$ and define numbers $n_1, \ldots, n_k$ by using the following constraints

$$1 < \lambda_{n_1}, \lambda_{n_2}^{-1} < (1/4)\delta\lambda_{n_1}^{-1}, \lambda_{n_3}^{-1} < (1/4)\delta\lambda_{n_2}^{-1}, \ldots, \lambda_{n_k}^{-1} < (1/4)\delta\lambda_{n_{k-1}}^{-1}, \delta < 0.1. \quad (3)$$

Now we take any cube with the edge $\delta$ and with the centre in some point $(\alpha_m)_{1 \le m \le k}$. Then the point $(\{t\lambda_{n_m}\})$ belongs to this cube, if

$$|\{t\lambda_{n_m}\} - \alpha_m| \le \frac{\delta}{2}. \quad (4)$$

Since the interval $(\alpha_m - \delta/2, \alpha_m + \delta/2)$ for sufficiently small $\delta$ has a length $< 0.1$ then the real numbers $t\lambda_{n_m}$ fractional parts of which lie in this interval have one and the same integral parts during continuous variation of $t$. So at $m = 1$, for some whole $r$, one has:

$$\frac{r + \alpha_1 - \delta/2}{\lambda_{n_1}} \le t \le \frac{r + \alpha_1 + \delta/2}{\lambda_{n_1}}. \quad (5)$$

The measure of a set of such values of $t$ does not exceed the size $\delta\lambda_{n_1}^{-1}$. The number of such intervals corresponding to different values of $r = [t\lambda_{n_1}] \le \lambda_{n_1}$ does not exceed

$$[\lambda_{n_1}] + 2 \le \lambda_{n_1} + 2.$$

So, the total measure of intervals satisfying (4) at $m = 1$ is less or equal to

$$(\lambda_{n_1} + 2)\delta\lambda_{n_1}^{-1} \le (1 + 2\lambda_{n_1}^{-1})\delta.$$

Consider now the case $m = 2$. Taking one of intervals of a view (4), we will have

$$\frac{s + \alpha_2 - \delta/2}{\lambda_{n_2}} \le t \le \frac{s + \alpha_2 + \delta/2}{\lambda_{n_2}}, \quad (6)$$

with some $s = [t\lambda_{n_2}] \le \lambda_{n_2}$. As we consider the conditions (4) for the values $m = 1$ and $m = 2$ simultaneously, we should estimate a total measure of intervals (6) which have nonempty intersections with intervals of a kind (5), using conditions (3). Every interval of a kind (6) is plased only in the one interval with the length $\lambda_{n_2}^{-1}$ where the expression $t\lambda_{n_2}$ has one and the same integral part $s$. The number of intervals with the length $\lambda_{n_2}^{-1}$ having a nonempty intersection with one fixed interval of a kind (5) does not exceed the size

$$[\delta\lambda_{n_1}^{-1}\lambda_{n_2}]+2 \le \delta\lambda_{n_1}^{-1}\lambda_{n_2}+2.$$

So, the measure of a set of values of $t$ for which intervals (6) have nonempty intersections only with one of intervals of a kind (5) is bounded by the value $(2+\delta\lambda_{n_1}^{-1}\lambda_{n_2})\delta\lambda_{n_2}^{-1}$. Since, the number of intervals (5) is no more than $\lambda_{n_1}+2$ then the measure of a set of values $t$ for which the condition (4) are satisfied simultaneously for $m=1$ and $m=2$ will be less or equal than

$$(\lambda_{n_1}+2)(2+\delta\lambda_{n_1}^{-1}\lambda_{n_2})\delta\lambda_{n_2}^{-1}.$$

It is possible to continue these reasoning considering all of conditions of a kind

$$\frac{l+\alpha-\delta/2}{\lambda_{n_m}} \le t \le \frac{l+\alpha+\delta/2}{\lambda_{n_m}}, m=1,\ldots,k.$$

Then we find the following estimation for the measure $m(\delta)$ of a set of those $t$ for which the points $(\{t\lambda_{n_m}\})$ located in the given cube with the edge $\delta$:

$$m(\delta) \le (2+\lambda_{n_1})(2+\delta\lambda_{n_1}^{-1}\lambda_{n_2})\cdots(2+\delta\lambda_{n_{k-1}}^{-1}\lambda_{n_k})\delta\lambda_{n_k}^{-1} \le \delta^k \prod_{m=1}^{\infty}(1+2m^{-2}).$$

Summarising over all such cubes, we receive the final estimation of a kind $\le 3c\varepsilon_1$ for the measure of a set of those $t$ for which $(\{t\lambda_{n_m}\}) \in B_1$ with the absolute constant

$$c = \prod_{m=1}^{\infty}(1+2m^{-2}).$$

We notice that the sequence $\Lambda = (\lambda_n)$, satisfying the conditions (3) defined above, depends on $\delta$. For each ball $B_k$ we fix some sequence $\Lambda_k$ using conditions (3). Considering all such balls, we designate $\Delta_0 = \{\Lambda_k \mid k=1,2,\ldots\}$.

Let $t \in E_0$ be any point. Taking $a = \{t\Lambda\}$ and applying the lemma 1 we conclude that the set $\Sigma(\{t\Lambda\})$ is contained in the finite union $\bigcup_{k \le n} B_k$ for some $n$ and

$$\Sigma(\{t\Lambda\}) \cap \left(B \setminus \bigcup_{k \le n} B_k\right) = \emptyset.$$

From here it follows that the set $E_0$ can be represented as a union of subsets $E_k, k=1,2,\ldots$, where

$$E_k = \{t \in E_0 \mid \Sigma(\{t\Lambda\}) \subset \bigcup_{s \le k} B_s\}.$$

So,

$$E_0 = \bigcup_{k=1}^{\infty} E_k; \quad E_k \subset E_{k+1} (k \geq 1).$$

Further, $m(E_0) = \lim_{k \to \infty} m(E_k)$, in agree with [4, p. 368]. As the set $E_k$ is a finite symmetric, then the measure of a set of values $t$, interesting us, is possible to estimate by using of any sequence $\Lambda_k$, since, as it has been shown above, the sets $\Sigma(\{t\Lambda\})$ for different values of $t$ have empty intersection. So,

$$m(E_k) \leq \limsup_{\Lambda' \in \Lambda_0} m(E_k(\Lambda')),$$

where $E_k(\Lambda') = \{t \in E_k \mid (\{t\Lambda'\}) \in \bigcup_{s \leq k} B_s\}$. Hence,

$$m(E_k(\Lambda')) \leq \sum_{s \leq k} m(E^{(s)}(\Lambda')),$$

where $E^{(s)}(\Lambda') = \{t \in E_0 \mid (\{t\Lambda'\}) \in B_s\}$. Applying the inequality found above, we receive:

$$m(E_k(\Lambda')) \leq 6c(\varepsilon_1 + \cdots + \varepsilon_k).$$

This result invariable for all $\Lambda' = \Lambda_r$ beginnig from some natural $r = r(k)$. Taking limsup, as $k \to \infty$, we receive the demanded result. The proof of the lemma 1 is finished.

### 4. Main theorem

The projection of the curve $(\{t\lambda_n\})_{n \geq 1}$ in two dimensional plane, i. e. the curve $(\{t\lambda_1\}, \{t\lambda_2\})$ has zero measure. By the theorem of Fubini, then, the product Lebesgue measure of the curve $(\{t\lambda_n\})_{n \geq 1}$ is also equal to zero. The following theorem shows that in the meaning of inroduced measure the behaviour of the curve $(\{t\lambda_n\})_{n \geq 1}$ is currently different.

**Theorem.** *Let the sequence $(\lambda_n)$ be an unbounded sequence of positive real numbers every finite subfamily of elements of which is linearly independent over the field of rational numbers. Then the curve $(\{t\lambda_n\}), t \in [0,1]$ is not $\mu_0$-measurable set in $\Omega$.*

*Proof.* Let's suppose the converse statement. Let the curve $(\{t\lambda_n\}), t \in [0,1]$ be measurable. Then it has a zero measure. Therefore, the union $U = \bigcup_{0 \leq t \leq 1} \Sigma(\{t\Lambda\})$ as a set constructed from the curve $(\{t\lambda_n\}), t \in [0,1]$ by an action of the group $\Sigma$ of all finite permutations has zero measure also, since it is a denumerable union of sets of zero measure. The set $U$ is a finite-symmetric. Let $n$ be any natural number. If we take a projection of the set $U$ to $\Omega$ by omitting the first $n$ coordinates (restricting the sequence $(\{t\lambda_n\})$), we get again the set $U_n$ of zero measure. Really, the set $U$ can be overlapped by the union of balls with the total measure of not exceeding $\varepsilon > 0$. Restricting the ball $B(\theta_0, \lambda)$ by omitting the first $n$ coordinates, and denoting the projection by

$S_N$, we get

$$S_N = \left\{(\theta_n) \mid \sum_{n=N+1}^{\infty} |\theta_n - \theta_n^0| e^{1-n} < \lambda \right\}.$$

Since

$$\sum_{n=N+1}^{\infty} |\theta_n - \theta_n^0| e^{1-N} = e^{-N} \sum_{n=1}^{\infty} |\theta_n - \theta_n^0| e^{1-n},$$

then denoting the projection of the point $\theta_0$ by $\theta_0'$, we have $S_N = B(\theta_0', e^N \lambda)$, and $e^N \lambda \to 0$ as $\lambda \to 0$ for any fixed $N$. From this one deduces the demanded statement.

Consider the sequence of sets $V_n = [0,1]^n \times U_n$, $n \in N$. It is obviously that $V_n \subset V_{n+1}$. Let $V = \bigcup_{n=1}^{\infty} V_n$. We have $\mu(V_n) = 0$ for all values of $n$. Therefore, also $\mu(V) = 0$, and the set $V$ is finite symmetrical. Then, there will be found some enumarable family of balls $B_r$ with a total measure not exseeding $\varepsilon$ the union of which contains the set $V$. For every fixed natural $n$ we define the set $\Sigma_n'(t\Lambda)$ as a closed set of all limit points of the sequence

$$\Sigma_n(\overline{\omega}) = \{\sigma\overline{\omega} \mid \sigma \in \Sigma \wedge \sigma(1) = 1 \wedge \cdots \wedge \sigma(n) = n\}.$$

Let

$$B^{(n)} = \{t \mid \{t\Lambda\} \in V \wedge \Sigma_n'(\{t\Lambda\}) \subset \bigcup_{r=1}^{\infty} B_r\}, \quad n = 1, 2, \ldots.$$

For every $t$ the sequence $\Sigma_{n+1}(\{t\Lambda\})$ is a subsequence of the sequence $\Sigma_n(\{t\Lambda\})$. Therefore, $\Sigma_{n+1}'(\{t\Lambda\}) \subset \Sigma_n'(\{t\Lambda\})$ and we have $B^{(n)} \subset B^{(n+1)}$. Then, one gets the inequality $m(B) \leq \sup_n m(B^{(n)})$ denoting $B = \bigcup_n B^{(n)}$.

Let's estimate $m(B^{(n)})$. The set $\Sigma_n'(\{t\Lambda\})$ is a closed set. Clearly, if we will "truncate" the sequences $\{t\Lambda\}$, remaining only components $\{t\lambda_k\}$ with indexes greater than $n$, and will denote the truncated sequence as $\{t\Lambda\}' \in \Omega$, then the set $\Sigma'(\{t\Lambda\}')$ also will be closed. Now we consider the products $[0,1]^n \times \{\{t\Lambda\}'\}$ (external brackets designate the set of one element) for every $t$. We have

$$\{t\Lambda\} \in [0,1]^n \times \{\{t\Lambda\}'\} \subset V.$$

Let $(\theta_1, \ldots, \theta_n) \in [0,1]^n$ is any point. There exist a neighborhood $V' \subset [0,1]^n$ of this point such that $(\theta_1, \ldots, \theta_n, \{t\Lambda\}') \in V' \times W \subset \bigcup_r B_r$, for some neighborhood $W$ of the point $\{t\Lambda\}'$. We, therefore, supplied every point $(\theta_1, \ldots, \theta_n) \in [0,1]^n$ with some pare of open sets $(V', W)$. Since the set $[0,1]^n$ is closed, then they can be found a finite number of open sets $V'$ the union of which contains $[0,1]^n$. The intersection of corresponding open sets $W$, being an open set, contains the point $\{t\Lambda\}'$. Therefore, for some finite set of indexes $R$ we have

$$[0,1]^n \times \{\{t\Lambda\}'\} \subset \bigcup V \times \bigcap W = [0,1]^n \times \bigcap W \subset \bigcup_{r \in R} B_r, \quad (7)$$

for each considered point $t$. It is clear that the set $R$ depends on the point $t$ and $\bigcap W \subset \bigcap_{r \in R} B'_r$ when $B'_r$ denotes the open set of trancated elements of $B_r$. The similar to (7) relationship is fair in the case when the point $\{t\Lambda\}$ would be replaced by any limit point $\overline{\omega}$ of the sequence $\Sigma(\{t\Lambda\})$ also, because $\overline{\omega} \in B_r$. If one denotes by $B'$ the union of all open sets of a kind $\bigcap_{r \in R} B'_r$ corresponding to every possible values of $t$ and of a limit point $\overline{\omega}$, we shall receive the relation

$$\{t\Lambda\} \in [0,1]^n \times \{\{t\Lambda\}'\} \subset A \subset [0,1]^n \times B' \subset \bigcup_{r=1}^{\infty} B_r,$$

for each considered values of $t$ and

$$\{\overline{\omega}\} \in [0,1]^n \times \{\overline{\omega}\}' \subset A \subset [0,1]^n \times B' \subset \bigcup_{r=1}^{\infty} B_r,$$

for each limit point $\overline{\omega}$. From this it follows the inequality

$$\mu_0^*\left([0,1]^n \times B'\right) = \mu_0^*(B') \leq \varepsilon,$$

where $\mu^*$ means an outer measure. The set $B'$ is open and $\Sigma'(\{t\Lambda\}') \in B'$. Now we can apply the lemma 1 and receive an estimation $m(B^{(n)}) \leq 6c\varepsilon$. Thus, we have $m(B) \leq 6c\varepsilon$. Since $\varepsilon$ could be chosen arbitrarily small then there exist $t$ such that $t \notin B$. So, $t \notin B^{(k)}$ for every $k = 1,2,...$. Consequently, for every $k$, there is a limit point $\overline{\omega}_k \in \Omega \setminus \bigcup_r B_r$ of the sequence $\Sigma_n(\{t\Lambda\})$. As the set $\Omega \setminus \bigcup_r B_r$ is closed, the limit point $\overline{\omega} = (\{t\Lambda\})$ of the sequence $(\overline{\omega}_k)$ will belong to the set $\Omega \setminus \bigcup_r B_r$. Therefore, $\{t\Lambda\} \notin \bigcup_{r \geq 1} B_r$ which contradicts the conditions of the theorem. Then the curve $(\{t\lambda_n\}), t \in [0,1]$ could not be $\mu_0$-measurable. The proof of the theorem is finished.

Mathematika 01/2009; 1(1).

8. Dzhabbarov I. Sh. On an identiy of harmonic analysis and its applications. Dokl. AS USSR, 1990, v.314, №5, 1052-1054.

9. R.Courant. Differential and integral calculus. M. Nauka, 1967. (rus).